\documentclass[10pt]{amsart}

 \usepackage{amssymb,latexsym,amsmath,amsthm}
\usepackage{graphicx} % support the \includegraphics command and options

\usepackage[numbers,sort&compress]{natbib}

\usepackage{fullpage}

 \renewcommand{\div}{\mathop{\mathrm{div}}\nolimits}
 \usepackage[numbers,sort&compress]{natbib}

\newtheorem*{thm*}{Theorem A}

\newtheorem{thm}{Theorem}[section]

\newtheorem{dfn}{Definition}[section]

\newtheorem{exam}{Example}[section]

\newtheorem{notation}{Notation}[section]

\newtheorem{lemma}{Lemma}[section]

\newtheorem{remark}{Remark}[section]

\numberwithin{equation}{section}

\usepackage{subcaption}
\usepackage{float}
\usepackage{fancyhdr} % This should be set AFTER setting up the page geometry
\begin{document}

\author{Mostafa Fazly}

\address{Department of Mathematics, The University of Texas at San Antonio, San Antonio, TX 78249, USA}
\email{mostafa.fazly@utsa.edu}

\author{Mark A. Lewis}
\address{Department of Mathematical  \& Statistical Sciences,  University of Alberta, Edmonton, AB T6G 2G1 Canada.}
\email{mark.lewis@ualberta.ca}

\author{Hao Wang}

\address{Department of Mathematical  \& Statistical Sciences,  University of Alberta, Edmonton, AB T6G 2G1 Canada.}
\email{hao8@ualberta.ca}

\thanks{MAL gratefully acknowledges a Canada Research Chair and an NSERC Discovery grant. HW gratefully acknowledge an NSERC Discovery grant. MF gratefully acknowledges a start-up grant of UTSA}

\def\IR{{\mathbb{R}}}

\newcommand{ \un }{\underline}

\newcommand{ \bg }{\begin{equation}}

\newcommand{ \ed }{\end{equation}}

% Omega    \Om
\newcommand{ \Om }{ \Omega}

% \partial Omega      \pOm
\newcommand{ \pOm}{\partial \Omega}

\title{Analysis of propagation for impulsive reaction-diffusion models}
\maketitle

\begin{abstract}
We study a hybrid impulsive reaction-advection-diffusion model given by a reaction-advection-diffusion equation composed with a discrete-time map in space dimension $n\in\mathbb N$.  The reaction-advection-diffusion equation takes the form
\begin{equation*}\label{}
u^{(m)}_t  = \div(A\nabla u^{(m)}-q u^{(m)})   + f(u^{(m)})  \quad  \text{for} \ \  (x,t)\in\mathbb R^n \times (0,1] ,
\end{equation*}
 for some function $f$,  a drift $q$ and  a diffusion matrix $A$. When the discrete-time map is local in space we use  $N_m(x)$ 
 to denote the density of population at a point $x$ at the beginning of reproductive season in the $m$th year and when the map is nonlocal we use $u_m(x)$.  The local discrete-time map is
\begin{eqnarray*}\label{}
\left\{ \begin{array}{lcl}
 u^{(m)}(x,0) = g(N_m(x))  \quad  \text{for} \ \ x\in \mathbb R^n ,
\\  N_{m+1}(x):=u^{(m)}(x,1)   \quad  \text{for} \ \ x\in \mathbb R^n ,
\end{array}\right.
 \end{eqnarray*}
for some function $g$.   The nonlocal discrete time map is
 \begin{eqnarray*}\label{}
\left\{ \begin{array}{lcl}
 u^{(m)}(x,0) = u_{m}(x) \quad  \text{for} \ \ x\in \mathbb R^n   , \\
\label{mainb2}  u_{m+1}(x) :=  g\left(\int_{\mathbb R^n} K(x-y)u^{(m)}(y,1) dy\right)   \quad  \text{for} \ \ x\in \mathbb R^n,
\end{array}\right.
 \end{eqnarray*}
 when $K$ is  a nonnegative normalized kernel.  Here, we analyze the above model from a variety of perspectives so as to understand the phenomenon of 
  propagation. We provide explicit formulae for the spreading speed of propagation in any direction $e\in\mathbb R^n$. Due to the structure of the model, we apply a simultaneous analysis of the differential equation and the recurrence relation to establish the existence of traveling wave solutions.  The remarkable point is that the roots of spreading speed formulas, as a function of drift, are exactly the values that yield blow-up for the critical domain dimensions, just as with the classical Fisher's equation with advection. We provide applications of our main results to impulsive  reaction-advection-diffusion models describing   periodically reproducing populations subject to climate change, insect populations in a stream environment with yearly reproduction,  and grass growing logistically in the savannah with asymmetric seed dispersal and  impacted by periodic fires.

 \end{abstract}


%\tableofcontents

\section{Introduction}

The hybrid models of a discrete-map and a reaction-advection-diffusion equation for species with
impulse and dispersal stages were proposed and studied in \cite{LL,flw}. 
Such discrete- and continuous-time hybrid models can describe  a seasonal event, such as reproduction or harvesting,  plus nonlinear turnover and dispersal throughout the year. We examine the propagation phenomenon for a general hybrid model  of the form given in \cite{flw}:
  \begin{equation}\label{main}
   u^{(m)}_t = \div(A\nabla u^{(m)}-q u^{(m)})  + f(u^{(m)})  \quad  \text{for} \ \  (x,t)\in\mathbb R^n\times (0,1], \ \ m=0,1,2,\cdots,
  \end{equation}
  where $A$ is a constant symmetric positive definite matrix and $ q$ is a constant vector.  In addition, we consider initial values either of the form
  \begin{eqnarray}\label{mainc}
 u^{(m)}(x,0) &=& g(N_m(x))  \quad  \text{for} \ \ x\in \mathbb R^n   , \\
 \label{mainc2} N_{m+1}(x) &:=& u^{(m)}(x,1)   \quad  \text{for} \ \ x\in \mathbb R^n,
  \end{eqnarray}
  or of the convolution form
\begin{eqnarray}\label{mainb}
 u^{(m)}(x,0) &=& u_{m}(x) \quad  \text{for} \ \ x\in \mathbb R^n   , \\
\label{mainb2}  u_{m+1}(x) &:= & g\left(\int_{\mathbb R^n} K(x-y)u^{(m)}(y,1) dy\right)   \quad  \text{for} \ \ x\in \mathbb R^n,
  \end{eqnarray}
  when $K$ is a nonnegative kernel, i.e. $K\ge 0$, and it is normalized so that is   $\int_{\mathbb R^n} K(x) dx =1$.  We do not
require $K$ to be symmetric so as to accommodate various applications. The Cauchy problem is given by (\ref{main})-(\ref{mainc2}) with $N_0(x)$ specified or by (\ref{main}), (\ref{mainb})-(\ref{mainb2}) with $u_0(x)$ specified. In \cite{flw}, authors studied the impulsive
reaction-diffusion model (\ref{main})-(\ref{mainc2}) and provided formulas for the critical domain sizes, associated with (\ref{Lstar}),
for various type domains. In the current article, we provide explicit formulas for the spreading speed of
propagation, which are counterparts of (\ref{cstara}), for such models.
 \begin{notation}
 Throughout this paper the matrix $A$ is defined as $A=(a_{i,j})_{i,j=1}^n$, the matrix $I=(\delta_{i,j})_{i,j=1}^n$ stands for the identity matrix,  and the vector $ q$ is $ q=(q_i)_{i=1}^n$. The matrix $A$ and the vector $ q$ have constant components unless otherwise is stated.   % We shall refer to vector fields $q$ with $\div q=0$ as divergence free vector fields.
 \end{notation}
%Throughout this article we make various assumptions on $f$ and $g$.  Occasionally, we
 We require the continuous function $g$ to satisfy the following assumption:
 \begin{enumerate}
 \item[(G0)]  $g$ is a positive function in $\mathbb R^+$,  $g(0)=0$,  $g'(0)>0$ and  $g(s)$ is nondecreasing in $0<s\le s^*$ for some $s^*>0$.  The quotient $g(s)/s$ is nonincreasing in $s>0$.  
  \end{enumerate}
  We now consider a few standard functions that satisfy the above assumptions.  Note that  the linear function
\begin{equation}\label{lg}
g(s)=\alpha s,
\end{equation}
 where $\alpha$ is a positive constant satisfies (G0).  For the special case where $\alpha=1$ and $f(s)=s(1-s)$,  model (\ref{main}) simplifies to the classical Fisher's equation for the spatial spread of an advantageous gene introduced  by Fisher in \cite{fish} and Kolmogorov, Petrowsky, and Piscounoff (KPP) in \cite{kpp} in 1937. The spreading speed formula was first computed by KPP and Fisher  for the semilinear parabolic equation 
 \begin{equation}\label{KPPF}
u_t-d\Delta u=f(u) \ \ \text{when} \ \ (x,t)\in\mathbb R^n\times \mathbb R^+.
\end{equation}
They proved that under certain assumptions on $f$, now called KPP nonlinearities, there is a threshold value $c^*=2\sqrt{d f'(0)}$ such that there is no front for $c<c^*$ and for all $c\ge c^*$ there is a unique front up to translations and dilations in terms of space and time, see also Aronson and Weinberger \cite{aw1,aw2}. With respect to  the standard Fisher's equation with the drift in one dimension,  
\begin{equation}\label{kppeq}
u_t = d u_{xx}-q u_x+f(u) \ \ \text{for} \ \ (x,t)\in \Omega\times\mathbb R^+,
\end{equation}
the critical domain size for the persistence versus extinction is
\begin{equation}\label{Lstar}
L^*:= \frac{2\pi d}{\sqrt{4d f'(0)- q^2}}  ,
\end{equation}
when $\Omega=(0,L)$ and the speeds of propagation to the right and left are
\begin{equation}\label{cstara}
 c_{\pm}^*(q)=2\sqrt{d f'(0)} \pm q,
 \end{equation}
when $\Omega=\mathbb R$.     The remarkable point, made by Speirs and Gurney \cite{sg}, is that $c_{\pm}^*(q)$ is a linear function of $q$ and  $L^*$ as a function of $q$ blows up to infinity exactly at roots of $c_{\pm}^*(q)$. For more information regarding the minimal domain size and analysis of reaction-diffusion models and connections between persistence criteria and propagation speeds, we refer interested readers to  Lewis et al. \cite{vll,plnl,LL,jl}, Lutscher et al. \cite{lnp}, 
 Murray and  Sperb \cite{ms},  Pachepsky et al. \cite{plnl},  Speirs and Gurney \cite{sg} and references therein.   For a broader perspective on reaction-advection-diffusion models in ecology we also refer interested readers to books of Cantrell and Cosner \cite{cc},  Fife \cite{f}, Lewis et al. \cite{lhl},   Murray \cite{m1,m2} and Perthame \cite{per}.  The  authors in \cite{flw} studied the impulsive reaction-diffusion model (\ref{main})-(\ref{mainc}) and provided formulas for the critical domain sizes, associated to (\ref{Lstar}), for various type domains. In particular, for   $n$-hypercube with the length $\Lambda_1=\cdots=\Lambda_n=\Lambda>0$, the critical domain size is given by 
\begin{equation}\label{Lstarn}
\Lambda^*:=\left\{ \begin{array}{lcl}
   2\pi d \sqrt \frac{n}{ 4d[f'(0)+\ln(g'(0))] - | q|^2},  &  \ \  \text{if} \ \  4d[f'(0)+\ln(g'(0))] - | q|^2>0,  \\
   \infty,  & \ \   \text{if} \ \  4d[f'(0)+\ln(g'(0))] - | q|^2<0.
\end{array}\right.
\end{equation}
In this article, we provide explicit formulas for the spreading speed of propagation, counterparts of (\ref{cstara}), for such models.

 One can consider nonlinear functions for $g$ such as the Ricker function,  that is
\begin{equation}\label{rg}
g(s)=s e^{\beta(1-s)},
\end{equation} where $\beta$ is a positive constant.  For the optimal stocking rates for fisheries,  mathematical biologists often apply  is the Ricker model \cite{rick},  introduced in 1954  to study  salmon populations with scramble competition for spawning sites leading to overcompensatory dynamics. %population in the  Pacific northwest, where they spawn in river beds and can spawn on top of a previous site.
The Ricker function is nondecreasing for $0<s \le s^*=\frac{1}{\beta}$ and satisfies assumption (G0).     Note also that the  Beverton-Holt function
\begin{equation}\label{bg}
g(s) =\frac{(1+\lambda)s}{1+\lambda s},
\end{equation} 
 with positive constant $\lambda$,   is an increasing function and satisfies assumption (G0).    %Functions of this form also go by the names of Michaelis-Menten, Monod and Holing type II.
This model was introduced to understand the dynamics of compensatory competition in  %in the context of
 fisheries by Beverton and Holt \cite{bh} in 1957.     Another example is the Skellam function
\begin{equation}\label{sg}
g(s)=\alpha(1-e^{-\beta s}),
\end{equation}   where $\alpha$ is a positive constant and $\beta>1$. This function satisfies  assumption (G0) and  was introduced by Skellam in 1951 in \cite{s} to study population density for animals, such as birds, which have contest competition for nesting sites, which leads to compensatory competition dynamics. %This function is also sometimes called   Ivlev function.
Note that the Skellam function behaves similar to the  Beverton-Holt function and it is nondecreasing for any $s>0$.  We shall use these functions in the application section (Section \ref{app}).  We refer interested readers to \cite{st} for more functional forms with biological applications.  

 We now provide some assumptions on the continuous function $f$. We assume that
 \begin{enumerate}
 \item[(F0)] $f(s)=f'(0)s+f_1(s)$  \text{where}  $f'(0)\neq 0$, $f_1(0)=f'_1(0)=0$ \text{and} $f_1\le 0$. 
 \end{enumerate}
 Note that we do not have any assumption on the sign of $f'(0)$.    Note that $f(s)=\alpha s-\beta s^2$,  $f(s)=\alpha s$ for $\alpha,\beta \in\mathbb R$  satisfy  (F0).   The above equation (\ref{main}) with (\ref{mainc})-(\ref{mainc2}) and (\ref{mainb})-(\ref{mainb2})  defines recurrence relations for $N_m(x)$ and  $u_m(x)$, respectively,  as
  \begin{equation}\label{oq}
  N_{m+1}(x)=Q[N_m(x)]   \quad  \text{for} \ \ x\in \mathbb{R}^n,
  \end{equation}
  and
   \begin{equation}\label{op}
  u_{m+1}(x)=P[u_{m}(x)]   \quad  \text{for} \ \ x\in \mathbb{R}^n,
  \end{equation}
  where $m\ge 0$ and $P$ and $Q$ are operators that depend on both the reaction-advection-diffusion equation and the discrete-time map, and thus $A, a,f $ and $ g$.  Most of the results provided in this paper are valid in any dimensions.  However we shall focus on the case of $n\le 3$ for applications.

This article is structured as follows. We provide properties of the recursion equation and spreading speeds $u_{m+1} = \Gamma[u_m]$ where $\Gamma$ is an operator on a certain set of functions on the habitat (Section \ref{secprop}).   Methods and techniques provided in this section are  based on those established by   Weinberger in \cite{w82, w}.  We consider the impulsive reaction-advection-diffusion equations with both local and nonlocal conditions.  We then provide explicit formulae for the spreading speed of  traveling waves and establish the existence of traveling wave (Section \ref{secspeed} and Section \ref{secspeedn}).   We also provide  applications of the main results  to models for  periodically reproducing  populations subject to climate change, to stream-dwelling insects that reproduce years, and to  grass growing logistically in the savannah, impacted by periodic fires (Section \ref{app}). Finally, we provide discussion and proofs for our main results.

\section{Formulation of the problem} \label{secprop}

In this section we study properties of the recursion
\begin{equation}\label{gam}
u_{m+1} = \Gamma[u_m]
\end{equation}
where $\Gamma$ is an operator on a certain set of functions on the habitat and $u_m (x)$ represents the gene fraction or population density at time $m$ at the point $x$
of the habitat. Later in the applications and proofs we shall set $\Gamma=P$ and $\Gamma=Q$,  where operators $P$ and $Q$ were introduced earlier. We shall define a wave speed $c^*(e)$ as a  scalar-valued function of unit direction vector  $e\in\mathbb R^n$ corresponding
to any operator $\Gamma$ which satisfies the following hypotheses, introduced by Weinberger as Hypotheses (3.1) in \cite{w82}. Such definitions and analysis of spreading speeds are discussed and extended in \cite{llw, w, bhn1, bhn2, aw1, aw2} and are used vastly in the literature that it is not limited to \cite{cc, cc01, hps, lhl, LL, llw, vll, jl}.

In all our applications the dependent variable $u$ has only nonnegative values. A population size will normally vary between 0 and some large upper bound $\pi_+$, but this upper bound could conceivably be $\infty$.  We define $\mathcal A$ as the set of non-negative continuous functions on $\mathbb R^n$ that are bounded by $\pi_+$. %We shall assume that $\Gamma[\alpha]>\alpha$ in some interval $(\pi_0,\pi_1) $ with $\Gamma[\pi_0]=\pi_0$ and $\Gamma[\pi_1]=\pi_1$, where $0\le \pi_0<\pi_+\le\pi_+$.
 Let us define the translation operator
 \begin{equation}
T_y[u(x)]=u(x-y).
\end{equation}
Note that a constant function is clearly translation invariant; that is $T_y[a]=a$ for $a\in\mathbb R$.  In addition we assume that $T_y[\Gamma[a]]=\Gamma[a]$.  Weinberger  \cite{w82} provided a list of properties for the operators in regards to the spreading speed
theory. In addition to these assumptions it is required that $\Gamma[u]$ to behave continuously with
respect to changes in $u$, and
\begin{enumerate}
\item[(H1)] $\Gamma[u]\in\mathcal A$ for all $u\in\mathcal A$.
\item[(H2)] $\Gamma[T_y[u]]=T_y[\Gamma[u]]$ for all $u\in\mathcal F$ and $y\in\mathbb R^n$.
\item[(H3)] There are constants  $0\le \pi_0<\pi_1\le \pi_+$ such that $\Gamma[s]>s$ for $s\in(\pi_0,\pi_1)$, $\Gamma[\pi_0]=\pi_0$, $\Gamma[\pi_1]=\pi_1$, if $\pi_1<\infty$.
\item[(H4)] If $u\le v$ then $\Gamma[u]\le \Gamma[v]$.
\item[(H5)] $u_m\to u$ uniformly on each bounded domain in $\mathbb R^n$ implies that $\Gamma[u_m]\to\Gamma[u]$ pointwise.
\end{enumerate}
%In order to define $c^*$, we begin by choosing a function $q(s)$ of one real variable with the properties
Now, we define the sequence $a_m(c,e;s)$ by the recursion
\begin{equation}
a_{m+1} (c,e;s)= \max\left\{ \phi(s), \Gamma [a_m(c,e;x\cdot e+s+c)](0) \right\}
\end{equation}
where $a_0=\phi(s)$ and $\phi(s)$ is continuous and nonincreasing function with $\phi(s)=0$ for $s\ge 0$ and $\phi(-\infty)\in (\pi_0,\pi_1)$ where $0\le \pi_0<\pi_1\le\pi_+$. From this definition, we notice that $a_0\le a_1$ and performing induction arguments we conclude that $a_m\le a_{m+1} \le \pi_+$.  In addition, the sequence $a_m$  is nondecreasing in $m$, nonincreasing in $s$ and
$c$, and continuous in $c$, $e$  and $s$. Therefore,  the sequence $a_m$ is monotone convergent  to a limit function $a(c,e;s)$ that is again nonincreasing in $s$ and
$c$, and bounded by $a_2$.  Applying Lemma 5.2 in \cite{w82}, we conclude that
\begin{equation}
a(c,e;-\infty)=\pi_+.
\end{equation}
However, the value of $a(c,e;\infty)$ may or may not be $\pi_+$. As it is shown in \cite{w82}  $a(c,e;\infty)=\pi_+$ if and only if there is $m$ that $a(c,e;0)>\phi(-\infty)$. We now define the spreading speed of propagation for in $e=(e_1,\cdots,e_n)\in \mathbb R^n$ direction  as $c^*(e)$   by
\begin{equation}
c^*(e)=\sup\{c; a(c,e;\infty)=\pi_+\} .
\end{equation}
Note that for the case of $a(c,e;\infty)=\pi_+$ for all $c$ we set $c^*(e)=\infty$.  The following theorem gives a method for  bounding the spreading speed of propagation.

\begin{thm}[Weinberger \cite{w82}] \label{thmw1}
If $m(x, dx)$ is a bounded nonnegative measure on $\mathbb R^n$ with the property that for all continuous functions $u$ with $0\le u\le \pi_1$,
\begin{equation}
\Gamma[u](x) \le \int_{\mathbb R^n} u(x-y) m(y, dy),
\end{equation}
then
\begin{equation}
c^*(e)\le \inf_{s>0} \left\{ \frac{1}{s} \ln  \int_{\mathbb R^n} e^{s x\cdot e} m(x, dx) \right\}.
\end{equation}
\end{thm}
\begin{thm}[Weinberger \cite{w82}] \label{thmw2}
If $m(x, dx)$ is a bounded nonnegative measure on $\mathbb R^n$ with the property that
\begin{equation}
\int m(x, dx) >1,
\end{equation}
and that, for all continuous positive bounded functions $u$ with $0\le u\le \pi_1$,
\begin{equation}
\Gamma[u](x) \ge  \int_{\mathbb R^n} u(x-y) m(y, dy) ,
\end{equation}
then,
\begin{equation}
c^*(e)\ge \inf_{s>0} \left\{  \frac{1}{s} \ln  \int_{\mathbb R^n} e^{s x\cdot e} m(x, dx)  \right\}.
\end{equation}
\end{thm}

    \section{Spreading Speed Formula; Linear Dynamics}\label{secspeed}

In this section, we provide an explicit formula for spreading speed of propagation for  (\ref{main}) with (\ref{mainc}) and (\ref{mainb})  in $e=(e_1,\cdots,e_n)\in \mathbb R^n$ direction. In addition,  we compute the direction $e$ such that the roots of the spreading speed as a function of advection coincide with the values for which critical domain dimensions  tend to infinity.

Suppose that $L[\cdot]$ is the linearization of operator $Q[\cdot]$ about zero. Define $M_0(x):=N_0(x)$,  then $M_{m+1}(x)=L[M_m(x)]$ where
  \begin{equation}\label{mainl}
   u^{(m)}_t +q \cdot \nabla u^{(m)} = \div(A\nabla u^{(m)})  + f'(0) u^{(m)}  \quad  \text{for} \ \  (x,t)\in\mathbb R^n\times (0,1]   ,
  \end{equation}
  satisfying
  \begin{eqnarray}\label{maincl}
 u^{(m)}(x,0) &=& g'(0) M_m(x)  \quad  \text{for} \ \ x\in \mathbb R^n   , \\
\label{maincl2}  M_{m+1}(x) &=& u^{(m)}(x,1)   \quad  \text{for} \ \ x\in \mathbb R^n .
  \end{eqnarray}
    Similarly, let $\bar L[\cdot]$ be the linearization of operator $P[\cdot]$ about zero. Set $\bar M_0(x):=u_0(x)$ and  $\bar M_{m+1}(x)=\bar L[\bar M_m(x)]$ where $u^{(m)}$ satisfies the linear equation (\ref{mainl}) with the following initial values
   \begin{eqnarray}\label{mainbl}
 u^{(m)}(x,0) &=& \bar M_{m}(x) \quad  \text{for} \ \ x\in \mathbb R^n   , \\
 \label{mainbl2}  \bar M_{m+1}(x) &=& g'(0)  \int_{\mathbb R^n} K(x-y)u^{(m)}(y,1) dy    \quad  \text{for} \ \ x\in \mathbb R^n. 
  \end{eqnarray}

  \begin{lemma}\label{propF}
Let  $L[\cdot]$ and $\bar L[\cdot]$ be the linearization of operator $Q[\cdot]$ and $P[\cdot]$  about zero, respectively. Then, for any $v\in C(\mathbb R^n)$ we have
\begin{equation}\label{ml}
  L(v)(x) =\int_{\mathbb R^n}  v(x-y) m(y,dy) \ \ \text{and} \ \ \bar L(v)(x) =\int_{\mathbb R^n}  v(x-y) l(y,dy),
\end{equation}
where the measures  $m$ and $l$ are defined as
\begin{eqnarray}\label{measm}
m(y,dy) &:=& g'(0) e^{ f'(0)} (2\pi)^{-n}  \frac{\pi^{\frac{n}{2}}}{\sqrt{\det A}}  e^{-\frac{1}{4} <A^{-1} \left(y- q\right),\left(y- q \right)>},\\
\label{measl}
l(y,dy) &:=& g'(0) e^{ f'(0)} (2\pi)^{-n}  \frac{\pi^{\frac{n}{2}}}{\sqrt{\det A}} \int_{\mathbb R^n} K(z-y)  e^{-\frac{1}{4} <A^{-1} \left(z- q\right),\left(z- q \right)>} dz,
\end{eqnarray}
and $<A^{-1}\eta,\eta>$ stands for $\eta^TA^{-1}\eta$ for any $\eta \in \mathbb R^n$.

  \end{lemma}
The proof of Lemma \ref{propF} is given in Section \ref{secproofs}.  Note that $m(y,dy)$ is a bounded nonnegative measure, since $A$ is a positive definite matrix.

 \begin{lemma}\label{lemcond}
 Let $m$ and $l$ be the measures introduced by (\ref{measm}) and (\ref{measl}). Then,
%  On the other hand,  it is straightforward calculation to show that
\begin{equation}\label{formule}
 \int_{\mathbb R^n}   m(x,dx) = \int_{\mathbb R^n}  l(x,dx) = g'(0) e^{ f'(0)} .
  \end{equation}
\end{lemma}

 \noindent {\bf Proof}.  This is a  special case of Lemma \ref{lemfor} for $s=0$. Note that from the assumptions on the kernel we have $$k(0)=\int_{\mathbb R^n} K(x)dx=1. $$

\hfill $ \Box$

Note that the linear operators  $L[v]$ and $\bar L[v]$ are explicitly formulated by (\ref{ml}) and therefore behave continuously with respect to changes in $v$. Here,  we emphasize on a few properties of these operators and in particular we show that assumptions (H1)-(H5) hold.  The operator $L$ commutes with any translation, meaning
\begin{equation}
L[T_y[v]](x)=T_y[L[v]](x),
\end{equation}
when $T$ is the shift operator $T_y[v](x)=v(x-y)$.  In addition, the comparison principle holds for the  operator $L$,  meaning that for $v_1\le v_2$ we have $L[v_1](x)\le L[v_2](x)$,  due to properties of the integral operator.   Let us recall that the integral operator is continuous with respect to its integrand. This elementary  fact implies that when  $v_n\to v$ as $n\to \infty$ uniformly then $L[v_n](x)\to L[v](x)$.  In the following theorem,  we determine
the traveling wave speeds in $e$ direction for the linear system. This result is a direct consequence of Theorem  \ref{thmw1} and Theorem  \ref{thmw2}.

\begin{lemma}\label{thmsp1} Let  $L[\cdot]$ and $\bar L[\cdot]$ be the linear operators given by (\ref{ml}) where
\begin{equation}\label{gfp}
g'(0) e^{ f'(0)}>1.
\end{equation}
   Then, the traveling wave speeds in $e$ direction  for these operators are
\begin{eqnarray}\label{upperc}
 c^*(e) &=& \inf_{s>0} \left\{ \frac{1}{s} \ln \int_{\mathbb R^n} e^{s x\cdot e} m(x,dx)  \right\} ,  \\
\label{upperc1}  c^*(e) &=& \inf_{s>0} \left\{\frac{1}{s} \ln \int_{\mathbb R^n} e^{s x\cdot e} l(x,dx)  \right\} ,
 \end{eqnarray}
respectively, when $m$ and $l$are the measures introduced by (\ref{measm}) and (\ref{measl}).
 \end{lemma}
 
% The ideas and method that we apply here are strongly motivated by the ones  derived by Weinberger in \cite{w82,w} and applied in \cite{llw,LL,lhl,vll}.
 \noindent {\bf Proof}. From the assumption (\ref{gfp}) and Lemma \ref{lemcond} we conclude that
\begin{equation}\label{}
 \int_{\mathbb R^n}   m(x,dx) = \int_{\mathbb R^n}   l(x,dx) >1 .
  \end{equation}
Since the linear operators $L$ and $\bar L$, given by Lemma \ref{propF}, satisfy assumptions (H1)-(H5), from Theorem \ref{thmw1} and  Theorem \ref{thmw2} we conclude the desired result.

\hfill $ \Box$

We now generalize (\ref{lemcond}) and compute the integrals in the right-hand side of the equations in (\ref{upperc}) and (\ref{upperc1}).
\begin{lemma}\label{lemfor}
 Let $m$ and $l$ be the measures introduced by (\ref{measm}) and (\ref{measl}). Then, for any $\mu\in\mathbb R$,
\begin{eqnarray}\label{formul}
 \int_{\mathbb R^n} e^{s x\cdot e}  m(x,dx) &=& g'(0) e^{ f'(0)} e^{s e \cdot  q+s^2 <Ae,e>} , \\
 \int_{\mathbb R^n} e^{s x\cdot e}  l(x,dx) &=&  g'(0) e^{ f'(0)} k(s) e^{s e \cdot  q+s^2 <Ae,e>} ,
  \end{eqnarray}
  when $k(s):=\int_{\mathbb R^n} K(x) e^{-s x\cdot e} dx$.
  \end{lemma}

%  \begin{notation}
%The symbol $< \cdot >$ stands for the dot product.
 % \end{notation}

The proof of Lemma \ref{lemfor} is given in Section \ref{secproofs}. We are now ready to provide an explicit formula for the spreading speed of the linearized model $L$.
\begin{thm}\label{thmsp}
Let $g'(0) e^{ f'(0)}>1$.  There exists a spreading speed $c_L^*$ associated with (\ref{mainl})-(\ref{maincl2}) of the following form
\begin{equation}\label{cse}
 c^*_{L}(e):=2\sqrt{<Ae,e>} \sqrt{f'(0)+\ln(g'(0))} +e\cdot q,
 \end{equation}
 such that initial data which is nonzero on a bounded set eventually spreads at speed $c_L^*$ in $e$-direction.

\end{thm}

\noindent  {\bf Proof}. In order to establish the formula of the traveling wave speed $c_L^*(e)$ for the linear operator $L$, we minimize the following function for $s>0$, as it is given in Lemma \ref{thmsp1} that is
\begin{equation*}\label{}
 c_L^*(e) = \inf_{s>0} \left\{ \frac{1}{s} \ln \int_{\mathbb R^n} e^{s x\cdot e} m(x,dx) \right\},
 \end{equation*}
where the measure $m$ is given by
\begin{equation*}\label{}
m(y,dy) = g'(0) e^{ f'(0)} (2\pi)^{-n}  \frac{\pi^{\frac{n}{2}}}{\sqrt{\det A}}  e^{-\frac{1}{4} <A^{-1} \left(y- q\right),\left(y- q \right)>}.
 \end{equation*}
Define $W(s):=\frac{1}{s} \ln \int_{\mathbb R^n} e^{s x\cdot e} m(x,dx)$ for $s>0$. Note that from  Lemma  \ref{lemfor}, we conclude that
$$ W(s)=  \frac{\ln[g'(0) e^{ f'(0)} ]}{s} + e \cdot  q+s <Ae,e>.  $$
Therefore,
\begin{equation*}\label{}
 c_L^*(e) = \inf_{s>0} W(s).
 \end{equation*}
 It is straightforward to compute the minimizer of the function $W$ as
 $$s^*:=\sqrt{\frac{\ln[g'(0)]+ f'(0) }{<Ae,e>} }.$$
Therefore, $c_L^*(e)=W(s^*)$ and this completes the proof.

\hfill $ \Box$

As a particular case, consider the speed of propagation $ c_L^*(e)$ in $e$-direction  as a function of the nonzero advection $q$  when $A=d(\delta_{i,j})_{i,j=1}^n$.   For the direction given by the $e=-\frac{q}{|q|}$, a unit vector, from (\ref{cse}) we have
\begin{equation}\label{}
 c_L^*\left(-\frac{q}{|q|}\right)=2\sqrt{d} \sqrt{f'(0)+\ln(g'(0))} -|q|.
  \end{equation}
Note that $c_L^*\left(-\frac{q}{|q|}\right)$ vanishes exactly at $|q|=2\sqrt{d[f'(0)+\ln(g'(0))]}$,  for which the critical domain dimensions in (\ref{Lstarn}) tend to infinity.  The following figure  clarifies this relation.

%%
%Figure 1
%%

%Figure 1%
\begin{figure}[H]
        \centering
        \begin{subfigure}[b]{0.45\textwidth}
                \includegraphics[width=\textwidth]{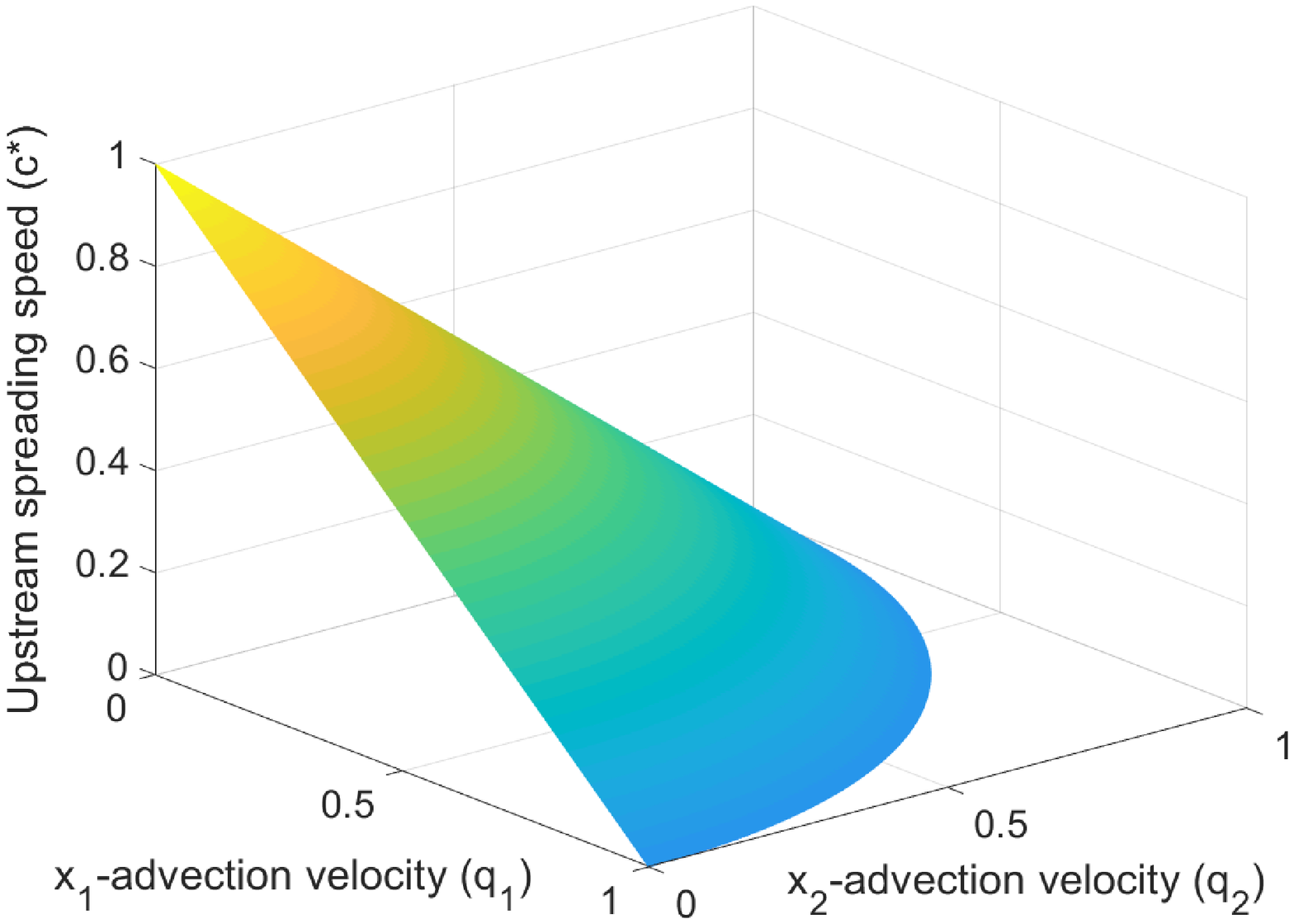}
                \caption{}
                \label{a}
        \end{subfigure}%
                        \begin{subfigure}[b]{0.45\textwidth}
                \includegraphics[width=\textwidth]{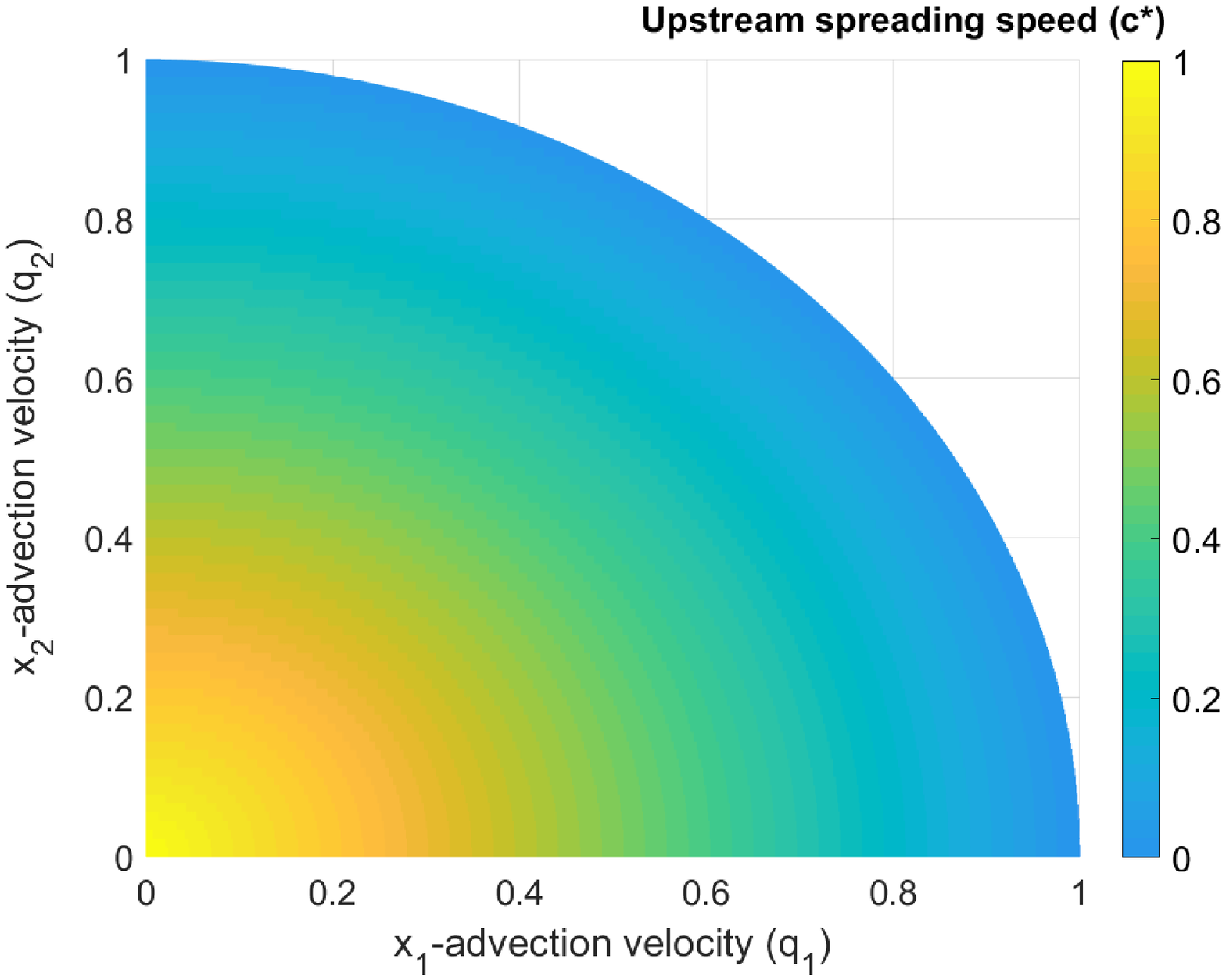}
                \caption{}
               % \label{b}
        \end{subfigure}
        \vspace*{-10pt}
        \caption{(A) Minimal speed $c_L^*(e)$ in the direction $e=-\frac{q}{|q|}$ when the advection velocity is given by $q=(q_1,q_2)$.  $c_L^*(e)$ is negative in the blank region.  The boundary between blank and shaded regions also gives the values of $q$ in (\ref{Lstarn}) where  the critical domain size $\Lambda^*$ becomes infinite. (B) Level sets of $c_L^*(e)$.} \label{stability}
\end{figure}

 \begin{remark} For the case of $n=1$,  $g(\tau)=\tau$ and $A=(d \delta_{ij})_{i,j=1}^n$, the spreading speed formula (\ref{cse}) turns into  the standard minimal speed
 \begin{equation}
 c_{\pm}^*=2\sqrt{d f'(0)} \pm q,
 \end{equation} in directions $e=\pm1$. The latter formula is the  upstream and downstream traveling wave speeds of the classical Fisher's equation.
  \end{remark}

The case with a Gaussian dispersal function $K$ also has a straightforward spreading speed formula. Suppose that $K(x)$ is a Gaussian distribution with mean $\mu$ and covariance matrix $B$ denoted by
\begin{equation}\label{KG}
K_G(x) : =\frac{1}{(4\pi)^{n/2} \sqrt{\det B}} e^{-\frac{1}{4} <B^{-1}(x-\mu),(x-\mu)>}.
\end{equation}
Here,  $\mu$ is a real vector and $B$ stands for some positive definite matrix.  We are now ready to provide an explicit formula for the spreading speed of the linearized model $\bar L$ associated with the above Gaussian distribution.

\begin{thm}\label{thmspk}
Let $g'(0) e^{ f'(0)}>1$. There exists a spreading speed $c_K^*$ associated to (\ref{mainl}) together with (\ref{mainbl})-(\ref{mainbl2}) of the following form
\begin{equation}\label{csk}
 c^*_{K_G}(e):=2\sqrt{<(A+B)e,e>} \sqrt{f'(0)+\ln(g'(0))} +e\cdot (q-\mu),
 \end{equation}
 such that initial data that is nonzero on a bounded set with positive measure eventually spreads at speed $c_K^*$ in $e$-direction. Here, $K_G$ is given by (\ref{KG}).  % and $e=(e_1,\cdots,e_n)$ is a unit vector in $\mathbb R^n$.
 \end{thm}
The proof of Theorem \ref{thmspk} is given in Section \ref{secproofs}.

    \section{Existence of Traveling Wave Solutions}\label{secspeedn}
    
     We start with the definition of traveling wave solutions.
  \begin{dfn} We say that $u_m(x)$ is a traveling wave solution of (\ref{gam}) in $e=(e_1,\cdots,e_n)\in \mathbb R^n$ direction if there exists a function $W$ and a constant $c$ such that $u_m(x)=W(x\cdot e- cm)$.
    \end{dfn}
    
  Suppose that $a=0$ and $A=0$ in (\ref{main}). Then $u(x,t)$ only depends on time and not on space, meaning that individuals do not advect or diffuse. Assume that $N_m$ represents the
number of individuals at the beginning of reproductive stage in the $m$th year. Then
 \begin{eqnarray}\label{ut}
 \left\{ \begin{array}{lcl}
    u_t (t)= f(u(t))  \quad  \text{for} \ \  t\in (0,1] , \\
 u(0)=g(N_m)    , \\
  N_{m+1}:=u(1).
\end{array}\right.
  \end{eqnarray}
   Separation of variables shows that
   \begin{equation}\label{equiseq}
  \int_{g(N_m)}^{N_{m+1}} \frac{d \omega}{f(\omega)}=1.
  \end{equation}
Note that a positive constant equilibrium of (\ref{ut})  satisfies
 \begin{equation}\label{equib}
 \int_{g(N)}^{N} \frac{d \omega}{f(\omega)}=1.
\end{equation}
Assume that $f$ satisfies (F0) and $g$ satisfies (G0) then
 \begin{equation}\label{Nstarine}
 1=\int_{g(N)}^{N} \frac{d \omega}{f(\omega)}\ge  \frac{1}{f'(0)}  \int_{g(N)}^{N} \frac{d \omega}{\omega}=  \frac{1}{f'(0)}  \ln \left |\frac{N}{g(N)} \right| \ge \frac{1}{f'(0)}  \ln \left |\frac{1}{g'(0)} \right|  .
   \end{equation}
In the light of above computations,   we assume that
 \begin{equation}\label{conditionfg}
 e^{f'(0)} g'(0)>1,
\end{equation}
 and an $N^*>0$ exists such that $f\neq 0$ on the closed interval with endpoints $N^*$ and $g(N^*)$ and
 \begin{equation}\label{nstar}
 \int_{g(N^*)}^{N^*} \frac{d \omega}{f(\omega)}=1.
\end{equation}

In what follows, we study the assumptions (H1)-(H5) for operators  $P$ and $Q$ given by (\ref{oq}) and  (\ref{op}). Let us start with the following order-preserving property.
\begin{lemma}
Assume that $v_0 $ and $ w_0$ are nonnengative continuous functions on $\mathbb R^n$ and $v_0(x)\le w_0(x)$ for all $x\in \mathbb R^n$.  Let $v(x,t)$ and $w(x,t)$ be solutions of (\ref{main}) with initial values of $v_0,w_0$ and  for operators $\Gamma=P$ and $\Gamma=Q$. Then, $v(x,t)\le w(x,t)$ for all $x\in\mathbb R^n$ and $t>0$.
%\begin{enumerate}
%\item $v(x,t)\le w(x,t)$ for all $x\in\mathbb R^n$ and $t>0$.
%\item $(K\star  v_0)(x) \le (K\star  w_0)(x) $ for all $x\in\mathbb R^n$.
%\item $g(v_0(x))\le g(w_0(x))$ for all $x\in\mathbb R^n$.
%\end{enumerate}
\end{lemma}
 The proof of this lemma follows from the comparison principle for reaction-diffusion equations, from the non-negativity of kernel $K$ and from the monotonicity
assumption on $g$. Therefore, (H1) and (H4) hold. It is straightforward to notice that operators $P$ and $Q$ commutes with all translations of the real line. For the convolution, one can apply the following change of variables argument,
\begin{equation}
\int_{\mathbb R^n} K(x-y-z) v(z) dz=\int_{\mathbb R^n} K(x-w) v(w-y) dw.
\end{equation}
This implies that (H2) holds for both $P$ and $Q$.   Note that for spatially constant solutions to (\ref{main}) when $u_0(x,0)=g(U_0)$ then the solution $u_m(x,0)$ remains spatially constant and satisfies
  \begin{equation}\label{main1}
  \frac{d}{dt} U_m =  f(U_{m})  \quad  \text{for} \ \  t\in (0,1]   ,
  \end{equation}
and $U_{m+1}(0)=g(U_m(1))$. Lewis and Li in \cite{LL} calculated the solution to this model explicitly for $f(u)=f'(0)u+f_1(u)$ when $f_1(u)=\gamma u^2$, and that is given by 
\begin{equation}
U_{m+1}= \frac{f'(0) g(U_m)}{(e^{-f'(0)}-1)\gamma g(U_m) + f'(0) e^{-f'(0)}}.
\end{equation}
 In the general case, one may refer to the computations as in (\ref{equiseq}) in this regard. 
Moreover, the following result was established for non-spatial model by Lewis and Li in \cite{LL} for $f_1(u)=\gamma u^2$ and by the authors in \cite{flw} for general nonlinearities $f$.
\begin{lemma}
\noindent 
\begin{enumerate}
\item If $g'(0) e^{f'(0)} \le 1$ and $U_0>0$, then $U_{m+1}\le U_m$ and $\lim_{m\to\infty} U_m=0$.
\item  If $g'(0)  e^{f'(0)}>1$ then there exists a unique $U^*>0$ with $Q(U^*)=U^*$.
\item  If $g'(0)  e^{f'(0)}>1$ and $0<U_0<U^*$, then $U_{m+1} > U_m$ and $\lim_{m\to\infty} U_m=U^*$.
\end{enumerate}
\end{lemma}
Applying similar arguments as in the above lemma one can see that
\begin{equation}
Q[0]=0, \ \ Q[U^*]=U^*, \ \ \text{and} \ \ Q[U]>U \ \ \text{for} \ \ U\in(0,U^*).
\end{equation}
Similar results hold for the operator $P$ and therefore assumption (H3) holds too. In order to show that (H5) is satisfied for operators $P$ and $Q$, we state that the integral operator
$$ \int_{\mathbb R^n} K(x-y) u(y) dy, $$
 with the kernel $K\in L^2(\mathbb R^n)$
is a compact operator in $\mathcal A\cap L^2(\mathbb R^n)$,  known also as the Hilbert-Schmidt operator. %Note that this operator is also known as the Hilbert-Schmidt operator with the kernel $K\in L^2(\mathbb R^n)$.
%ntegral operators with Hilbert-Schmidt kernel are compact
%\begin{lemma}\label{lemcon}
%The convolution operator $\int_{\mathbb R^n} K(x-y) u(y) dy$ is compact in  $\mathcal A$.
%\end{lemma}
Since  $g$ is continuous, this implies that both operators $P$ and $Q$ are  compact in $\mathcal A$ on every bounded set.  

We  now provide the existence of a traveling wave solution of speed $c$ in $e$ direction whenever $c \ge c^*(e)$. Here,  $c^*(e)$ is given by (\ref{upperc}) and (\ref{upperc1}) for the linearized operators $L[\cdot]$ and $\bar L[\cdot]$ stated in (\ref{mainl}) with conditions (\ref{maincl})-(\ref{maincl2}) and (\ref{mainbl})-(\ref{mainbl2}), respectively. We say that  the recursion model (\ref{main}) together with either (\ref{mainc})-(\ref{mainc2}) or (\ref{mainb})-(\ref{mainb2}) is linearly determinate in the direction of $e$, since  $c^*(e)$ is determined exactly in terms of the behavior of the linearized system. 

\begin{thm}\label{thmcw}
For every $c \ge c^*(e)$ there exists a traveling wave solution of speed $c$ that is $W(x\cdot e-mc)$  for model (\ref{main}) together with either (\ref{mainc})-(\ref{mainc2}) or (\ref{mainb})-(\ref{mainb2}). Here,   $W$ is a nonincreasing function that $W(-\infty)=\pi_1$ and $W(+\infty)=0$. In addition,  the system is linearly determinate; that is $c^*(e)$ is given by (\ref{upperc}) and (\ref{upperc1}).
\end{thm}

 %In one dimension that is $n=1$ for $g(N)=N$ and $A=d$, the formula (\ref{cse}) when $e=\pm1$ is the standard spreading speed $c_{\pm}^*$ for the Fisher's equation provided in (\ref{cstarass}).

%%%%%%%%%
%%%%%%%%%%%%
%%%%
As an example,  we  compute the spreading speed for particular models in two dimensions.
\begin{exam}\label{Exam1} 
 Let $n=2$, $e=(\cos\theta,\sin\theta)$, $A=(a^2_{ij}\delta_{ij})_{i,j=1}^n$ and $ q=(q_1,q_2)$ where $0\le \theta < 2 \pi$ and $a_{ij}$ and $q_i$ are constant. Applying formula (\ref{cse}) we get
\begin{equation}\label{}
c^*(e)=2\sqrt{a^2_{11} \cos^2 \theta +a^2_{22} \sin^2 \theta }  \sqrt{f'(0)+\ln(g'(0))} +q_1\cos \theta+ q_2 \sin\theta.
  \end{equation}
\end{exam}
This formula clarifies the dependence of the spreading speed of propagation, for any angle $\theta$,  on anisotropic diffusion coefficients $a_{11}$ and $a_{22}$ and on advection coefficients $q_1$ and $q_2$.    Define the following set
$$ \mathcal S:=\{x\in\mathbb R^n, c^*( \tilde  e)\ge \tilde e\cdot x \ \ \text{for all vectors} \ \tilde e\}.$$
Suppose that $c^*(e)$ is the propagation speed in the $e$-direction of a traveling wave, then the convex set $\mathcal S$ would be the ray surface.  As it is mentioned by Weinberger in \cite{w}, the ray speed in $ e$-direction is defined to be the largest value of $\alpha$ such that $\alpha  e \in \mathcal S$ and it is related to $c^*$ by the following formula
 \begin{equation}\label{Cemin}
C( e) := \min_{e\cdot \tilde e>0, \tilde e \in \mathbb R^n} \frac{ c^*( \tilde  e)}{e\cdot\tilde e}.
  \end{equation}
 In other words, the ray speed in direction $e$, $C( e), $ is the minimizer of the speed in the direction $ e$ among all the fronts, even those in directions $\tilde e\neq e$, providing $\tilde e \cdot e>0$. Thus $C(e)$ is the minimum possible rate of spread in the $e$-direction, given that an expansion front occurs in a direction whose projection onto $e$ is positive.  This ray speed will be less than or equal to the speed of an expansion front travelling in the $e$-direction $c^*(e)$. Note that the above formula is also known as Freidlin-G\"{a}rtner's formula \cite{gf}. Just as in the above example, we shall compute the ray speed in the direction $e=(\cos\theta,\sin\theta)$ as it follows. Let  $ q=0$ in Example \ref{Exam1}. Then, for a unit vector $\tilde e=(\tilde e_1,\tilde e_2)\in\mathbb R^2$, we have 
 \begin{equation}\label{}
c^*(\tilde e)=2\sqrt{a^2_{11} \tilde e^2_1 +a^2_{22} \tilde e^2_2 }  \sqrt{f'(0)+\ln(g'(0))}.
  \end{equation}
 From this and (\ref{Cemin}), we conclude 
  \begin{equation}\label{Ce1}
C(e)=  2\sqrt{f'(0)+\ln(g'(0))} \min_{\cos (\theta) \tilde e_1 + \sin (\theta) \tilde e_2 >0} \ \   \frac{ \sqrt{a^2_{11} \tilde e^2_1 +a^2_{22} \tilde e^2_2 }  }{\cos (\theta) \tilde e_1 + \sin (\theta) \tilde e_2 }.
  \end{equation}
It is straightforward to notice that for all units vector $\tilde e=(\tilde e_1,\tilde e_2)\in\mathbb R^2$ and for $0\le \theta\le \pi/2$, the minimum is attained and it is given by 
  \begin{equation}\label{Ce1m}
\min_{\cos (\theta) \tilde e_1 + \sin (\theta) \tilde e_2 >0} \frac{ \sqrt{a^2_{11} \tilde e^2_1 +a^2_{22} \tilde e^2_2 }  }{\cos (\theta) \tilde e_1 + \sin (\theta) \tilde e_2 } = \sqrt{\frac{a_{11}^2a_{22}^2}{ a^2_{11} \sin^2  \theta  +a^2_{22} \cos^2 \theta }  } . 
  \end{equation}
Therefore, 
\begin{equation}\label{Ce}
C(e)= 2 \sqrt{\frac{a_{11}^2a_{22}^2}{ a^2_{11} \sin^2 \theta +a^2_{22} \cos^2 \theta }  }   \sqrt{f'(0)+\ln(g'(0))} . 
\end{equation}
Note that both $c^*(e)$ and  $C(e)$  as  functions of $\theta$  for $0\le \theta\le \pi/2$ are either decreasing or increasing depending on $a_{11}>a_{22}$ or $a_{22}>a_{11}$. In addition,
$c^*(e)\ge C(e)$ if and only if 
$$(a^2_{11}-a^2_{22})^2\sin^2 \theta \cos^2 \theta \ge 0.$$
 Therefore, $C(e)=c^*(e)$ if and only if $e=(\pm 1,0)$ or $e=(0,\pm 1)$ or $a_{11}=a_{22}$. In the following figures we illustrate the relation between $c^*(e)$ and  $C(e)$.
%%%%%%%%%%%%%%%%%%%%
%%%%%%%%%%%%%%%%%
%%%Figure 2
\begin{figure}[H]        \centering   \begin{subfigure}[b]{0.42\textwidth}  \includegraphics[width=\textwidth]{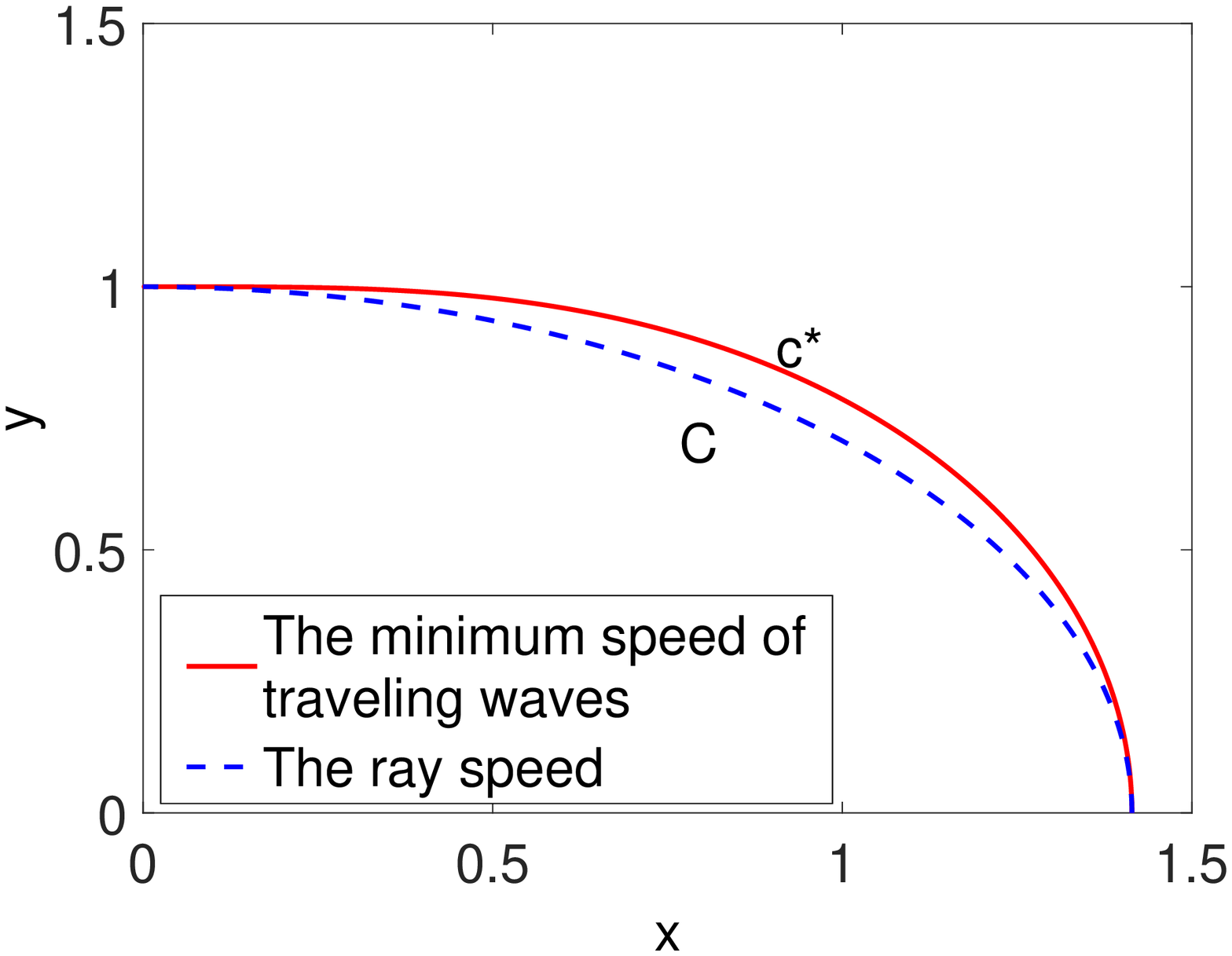} \caption{}  \label{a} \end{subfigure} \begin{subfigure}[b]{0.47\textwidth}  \includegraphics[width=\textwidth]{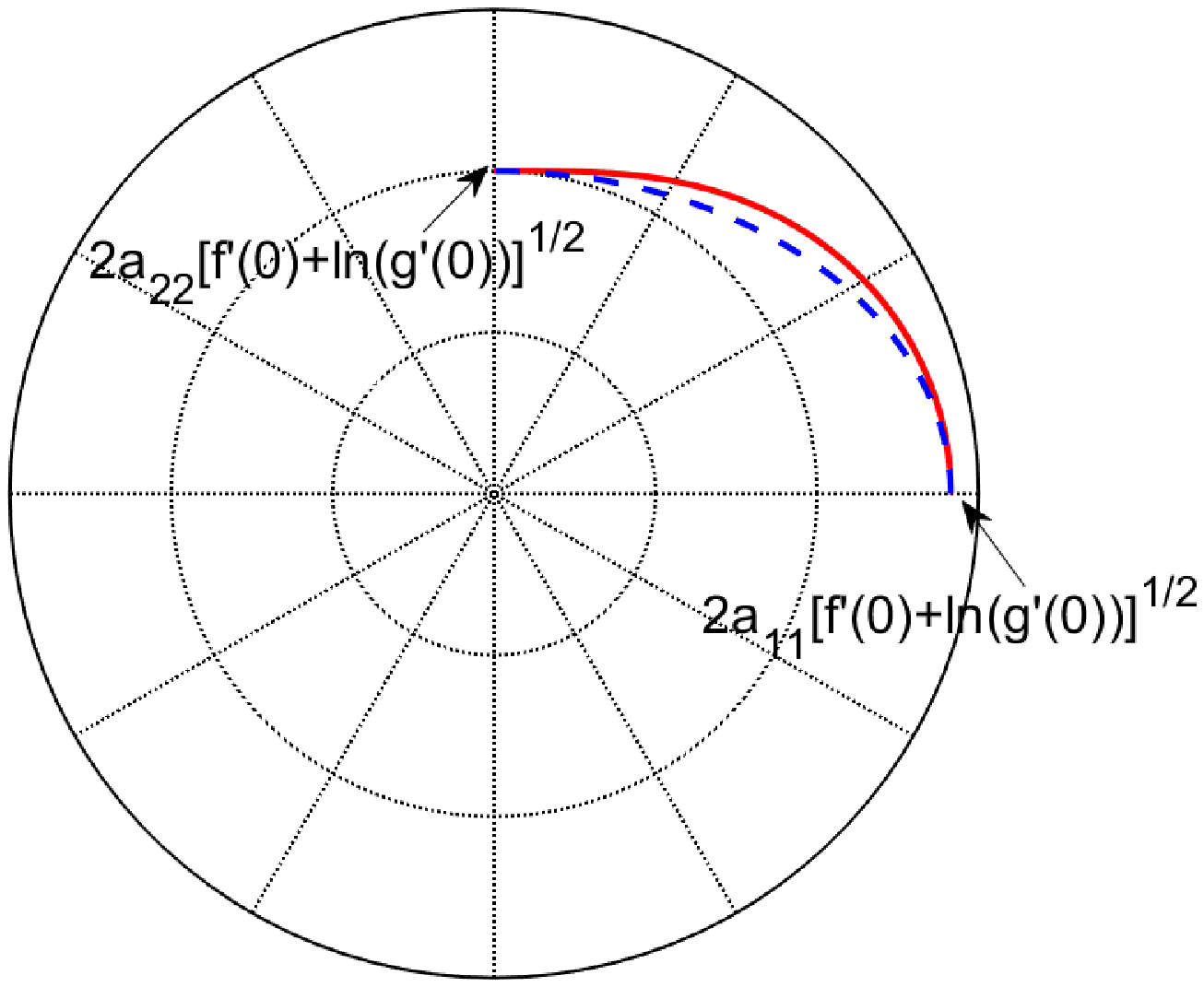} \caption{}        \end{subfigure} \vspace*{-10pt} \caption{Minimal speed and Ray speed in $e=(\cos \theta,\sin \theta)$ direction in $2D$ when $0\le \theta\le \frac{\pi}{2}$.  In (A) the left the case where $a_{11}>a_{22}$ is shown.  In (B) the quantitites $c^*$ (red) and $C$ (blue) are shown in polar coordinates. } \label{stability} \end{figure}
%\vspace*{-20pt}
%\begin{figure}[H]\begin{center} \includegraphics[height=1.8in]{Figure3.eps} \vspace*{-10pt} \caption[Indirect incidence term]{Minimal and Ray speeds in polar coordinates.} \label{fig:Figure1} \end{center} \end{figure}

%We end this section with providing asymptotic behavior of the traveling speed (\ref{csk}) for the cases of large drift, large diffusion and strong anisotropy. Note that $c^*(e)$ is a continuous function in terms of components of $A,q$ and $f'(0),g'(0)$.\begin{remark}\noindent \begin{enumerate}\item Large drift; let the advection be $q_M=M q$ where $q$ is a vector.   Then,$$ \lim_{M\to\infty} \frac{c^*(e)}{M}=e\cdot q.$$\item Large diffusion; let the diffusion matrix be $A_M=M^2 I$ where $I$ is the identity matrix. Then,$$ \lim_{M\to\infty} \frac{c^*(e)}{M}=2\sqrt{[f'(0)+\ln(g'(0))]}.$$\item Strong anisotropy; let the diffusion matrix be  $A=(a_{ij}\delta_{ij})_{i,j=1}^{n}$ and $a_{kk}\to 0$ for some index $k$. Then  $$\lim_{a_{kk}\to 0} c^*(e) = 2\sqrt{\Sigma_{i\neq k} a^2_{ii} e_i^2} \sqrt{f'(0)+\ln(g'(0))}.$$\end{enumerate} \end{remark}

\section{Applications}\label{app}

In this section, we provide applications of main theorems regarding analysis of spreading speed formula provided in Section \ref{secspeed} and Section \ref{secspeedn}.

\subsection{Population subject to climate change}
Climate change, especially global warming, has greatly changed the distribution and habitats of biological species. Uncovering the potential impact of climate change on biota is an important task for modelers \cite{pl}. We investigate the dependence of the critical domain size and the spreading speed on the climate shifting speed.

We consider a rectangular domain that is $\Omega=[0,L_1]\times[0,L_2]$ moving in the positive $x$-axis direction at speed $c$. Outside this domain conditions are hostile to population growth, while inside the domain there is random media, mortality and periodic reproduction. Using the approach of \cite{pl} this problem is transformed to a related problem on a stationary domain. Consider the model
\begin{eqnarray}\label{mainmr1}
 \left\{ \begin{array}{lcl}
    u^{(m)}_t = d \Delta u^{(m)}  -q\cdot \nabla u^{(m)} - \gamma u^{(m)} \quad  \text{for} \ \  (x,t)\in \Omega\times(0,1] ,\\
    u^{(m)}(x,t)=0   \quad  \text{for} \ \  (x,t)\in\partial\Omega\times (0,1], \\
 u^{(m)}(x,0)=g(N_m(x))  \quad  \text{for} \ \ x\in \Omega,\\
  N_{m+1}(x):=u^{(m)}(x,1)   \quad  \text{for} \ \ x\in \Omega,
\end{array}\right.
  \end{eqnarray}
  where $g$ is  the Beverton-Holt function that is  $g(N_m) =\frac{(1+\lambda)N_m}{1+\lambda N_m}$ and  $q=(c,0)$.  Similar to previous examples the assumption (\ref{gfp}) implies that  $\ln(1+\lambda)>\gamma$.   From the critical domain dimensions provided by the authors in \cite{flw} we conclude that
  \begin{equation}\label{L1L2}
  \frac{1}{L_1^2} + \frac{1}{L_2^2}<  \frac{1}{d\pi^2} \left[\ln(1+\lambda)-\gamma - \frac{c^2}{4d}\right],
  \end{equation}
when   $\ln(1+\lambda)>\gamma + \frac{c^2}{4d}$ we have  $\lim_{m\to \infty} N_m(x)=\bar N(x)$ for some positive equilibrium $\bar N(x)$ which refers to persistence of population.  In other words, (\ref{L1L2}) yields the parameter $c$ must be bounded by
\begin{equation}\label{cL}
  c^2 < 4d[\ln(1+\lambda)- \gamma ] - (2d\pi)^2 \left[\frac{L_1^2+L_2^2}{L_1^2L_2^2} \right] .%<4d[\ln(1+\lambda)- \gamma ] .
\end{equation}
  Let us mention that Theorem \ref{thmsp} implies that the speed of propagation for the model (\ref{mainmr1}) on an infinite domain $\Omega=\mathbb R$ in the $x-$axis direction, $e=1$, and in the opposite of $x-$axis direction, $e=-1$, is given by 
    \begin{equation}\label{cstar}
c^{*}=2 \sqrt{d [\ln(1+\lambda)- \gamma]} \pm c.
  \end{equation}
  Note that when $c$ satisfies (\ref{cL}) then the speed of propagation in (\ref{cstar}) is positive.  This implies that  persistence and ability
to propagate should be closely connected.  This is the case from a biological perspective as well. For example, if a population
cannot propagate upstream but is washed downstream,  it will not persist.   We refer interested readers to Speirs and  Gurney \cite{sg} and Pachepsky et al. \cite{plnl} for similar arguments regarding  Fisher's equation with advection. Note also that   the speed of propagation vanishes in (\ref{cstar})  exactly at the values that make the right-hand side of (\ref{L1L2}) zero.  One can compare this relation to the one given in (\ref{Lstar}) and (\ref{cstara}) for the classical Fisher's equation.

\subsection{Population of a stream insect species (e.g. stoneflies, mayflies)}
We now apply the nonlocal model to formulate a mixed continuous-discrete model for a single population of a stream
insect species (e.g. stoneflies, mayflies) with two distinct, non-overlapping developmental
stages. Such models are studied by Vasilyeva, Lutscher and Lewis  in \cite{vll} and  by Speirs and Gurney in \cite{sg}. Let $ u^{(m)}(x,t)$ be the density of the larval population at time $t$ and location $x$ during season
$m$. Larvae are transported by diffusion  rate $d$ and by drift
speed $q$, and experience possibly density-dependent death according to some
positive function $f$.  Then, consider the model
\begin{eqnarray}\label{mainmr}
 \left\{ \begin{array}{lcl}
     u^{(m)}_t = d \Delta  u^{(m)}  -q\cdot \nabla  u^{(m)} - r u^{(m)} \quad  \text{for} \ \  (x,t)\in\mathbb R\times (0,1]  ,\\
u^{(m)}(x,0) = u_{m}(x) \quad  \text{for} \ \ x\in \mathbb R   , \\
  u_{m+1}(x) =  g\left(\int_{\mathbb R^n} K(x-y)u^{(m)}(y,1) dy\right)   \quad  \text{for}\ \ x\in \mathbb R,
\end{array}\right.
  \end{eqnarray}
  where $r$ is the natural mortality rate, $d$ is the diffusion rate,  $g$ is  the Beverton-Holt function $g(\tau) =\frac{(1+\lambda) \tau }{1+\lambda \tau}$ for $\lambda>0$ and $q$ is the speed rate of drift.  In this model, the adult dispersal stage is modeled by a Gaussian normal distribution with mean  up-river dispersal distance  $\mu$ and covariance $\sigma^2$, denoted by (\ref{KG}) in one dimension, that is
  \begin{equation}\label{}
K(x) : =\frac{1}{ \sqrt{4\sigma^2 \pi} } e^{-\frac{1}{4\sigma^2} \left |x-\mu\right|^2}.
\end{equation}
  Note that the condition (\ref{gfp}) is equivalent to
  \begin{equation}\label{cond1l}
(1+\lambda) e^{-r}>1.
\end{equation}
Theorem \ref{thmspk} implies that
\begin{equation}\label{}
 c^*:=2\sqrt{d+\sigma^2} \sqrt{\ln(1+\lambda)-r} \pm  (q-\mu) .
 \end{equation}
Note that $c^*>0$ is equivalent to
\begin{equation}\label{}
\lambda+1>e^{ r+\frac{(q-\mu)^2}{4(d+\sigma^2)} }.
 \end{equation}
 Therefore, the population spreads in both directions exactly when the above condition is satisfied. This is, in general, stronger than the non-extinction condition (\ref{cond1l}) hat is $(1+\lambda) e^{-r}>1.$ When a population can persist on a bounded domain, it can spread upstream and
downstream on the infinite domain.

\subsection{Grass growing logistically in the savannah  with asymmetric seed dispersal and  impacted by periodic fires.}
The high grass productivity on savannahs can lead to more frequent savannah fires. We apply the  impulsive model (\ref{main})-(\ref{mainc}) with Fisher-KPP nonlinearity to study grass growing logistically in the savannah, impacted by periodic fires. For an isotropic impulsive model in one dimension, see \cite{yctbd}.  Seed dispersal is the movement or transport of seeds away from the parent plant  and it is often directionally biased, because of the inherent directionality of
wind and many other dispersal vectors. Therefore, we consider an anisotropic model to characterize seed dispersal with
advection in two dimensions.  Let $n=2$, $e=(\cos\theta,\sin\theta)$, $A=(a^2_{ij}\delta_{ij})_{i,j=1}^n$ when $a_{11}\neq a_{22}$ and $ q=(q_1,q_2)$ where $0\le \theta < 2 \pi$ and $a_{ij}$ and $q_i$ are constant. For $(x,y,t)\in\mathbb R\times \mathbb R \times (0,1]$, consider
\begin{eqnarray}\label{umta}
 \left\{ \begin{array}{lcl}
     u^{(m)}_t =  a_{11} u_{xx}^{(m)} + a_{22} u_{yy}^{(m)} - <q_1,q_2>\cdot <  u_x^{(m)},u_y^{(m)}> + r(1-u^{(m)})u^{(m)} ,\\
u^{(m)}(x,y,0) =g(N_m(x,y)) \quad  \text{for} \ \ x,y\in \mathbb R \times \mathbb R  , \\
  N_{m+1}(x,y)=u^{(m)}(x,y,1)  \quad  \text{for}\ x,y\in \mathbb R \times \mathbb R,
\end{array}\right.
  \end{eqnarray}
   where  $g(N_m) =(1-s) N_m$ and $0<s<1$.  The assumption (\ref{conditionfg}) is equivalent to
  \begin{equation}\label{ers}
  e^r (1-s) > 1.
  \end{equation}
  Assuming that (\ref{ers}) holds,  straightforward calculations show that the positive equilibrium $N^*$ that solves (\ref{nstar}) is of the form
  \begin{equation}\label{Nstarr}
  N^*:= \frac{ (1-s)e^r -1 }{(1-s) (e^r-1)}.
   \end{equation}
   From (\ref{cse}), we obtain an explicit formula of spreading speed of propagation in direction $e=(\cos\theta,\sin\theta)$ of the form of
\begin{equation}\label{}
c^*(e)=2\sqrt{a^2_{11} \cos^2 \theta +a^2_{22} \sin^2 \theta }  \sqrt{\ln (1-s)+r} +q_1\cos \theta+ q_2 \sin\theta.
  \end{equation}
We now compute the ray speed in the direction $e=(\cos\theta,\sin\theta)$ when $ q=0$. Then
\begin{equation}\label{}
C(e)= 2 \sqrt{\frac{a_{11}^2a_{22}^2}{ a^2_{11} \sin^2 \theta +a^2_{22} \cos^2 \theta }  }   \sqrt{\ln (1-s)+r} .   \end{equation}
Straightforward computations show that  $c^*(e)\ge C(e)$ if and only if $(a^2_{11}-a^2_{22})^2\sin^2\theta \cos^2 \theta\ge 0.$  Because of anisotropic diffusion, the above inequality implies the wave speed is strictly larger than  the ray speed in any direction that is not parallel to $x$-axis, that is $e=(\pm 1,0)$, or $y$-axis, that is $e=(0,\pm 1)$.
%Therefore, $C(e)=c^*(e)$ if and only if $e=(\pm 1,0)$ or $e=(0,\pm 1)$ or $a_{11}=a_{22}$.

\section{Discussion} \label{secdis}
We analyzed population spread and traveling waves for impulsive reaction-advection-diffusion equation models for species with distinct reproductive and dispersal stages on spatial domains $\Omega\subset \mathbb R^n$ when $ n\ge1$.  The case of one dimensional model was studied by Lewis and Li in \cite{LL} and Vasilyeva et al. in \cite{vll} and the critical domain size  for higher dimensions by the authors in \cite{flw}.   Unlike the analysis of standard partial differential equation models,  the study of impulsive reaction-advection-diffusion models requires a simultaneous analysis of the differential equation and the recurrence relation.  This fundamental fact rules out certain standard mathematical analysis theories for analyzing solutions of these type models,  but it  opens  up various new ways to apply the model.   The models can be variously considered as a description for a continuously growing and a dispersing population with pulse harvesting, a dispersing population with periodic reproduction, or  a population with individuals immobile during the winter.

On the entire space $\mathbb R^n$,  we provided an explicit formula for the spreading speed of propagation in any direction $e\in \mathbb R^n$ in terms of the same set of model parameters used for computing critical domain sizes and extreme volume sizes.  Our applications section demonstrates that the  possible modeling questions that can be addressed are broadranging, and in many cases need to be addressed in higher dimensions than dimension one.  Our paper has developed the tools to make this possible. Study of the minimal speed of propagation and the asymptotic spreading speed has attracted the attention of many mathematicians and scientists for the past few decades, see \cite{aw1,aw2,gf,w82,w,f}.    Many authors have driven formulae for the spreading speed  of propagation for parabolic equations with a non constant diffusion matrix  and a non constant advection vector field on cylinders,  periodic domains and general domains, see  G\"{a}rtner-Freidlin  \cite{gf}, Mallordy-Roquejoffre \cite{fr}, Heinze-Papanicolaou-Stevens \cite{hps}, Berestycki-Hamel-Nadirashvili  \cite{bhn1,bhn2} and references therein. In these articles, authors used  variational principles, more precisely min-max theories, to express the spreading speed formulae.  One can apply the ideas and mathematical techniques used in these references to develop a theory of  propagation for the impulsive reaction-diffusion equation models with a nonconstant diffusion matrix  and a nonconstant advection vector with hostile and flux boundary conditions.

\section{Proofs}\label{secproofs}

In this section, we provide proofs for the main results of Section \ref{secspeed} and Section \ref{secspeedn}. We shall start with the following technical lemma that is used frequently in the proofs.  We omit the proof of this lemma since it is straightforward.

 \begin{lemma}\label{lemcom}
Let  $A=(a_{i,j})_{i,j=1}^n$ be a symmetric positive definite matrix  with constant components.  Then,  the following formulas hold: 
\begin{equation}\label{intA}
\int_{\mathbb R^n} e^{i z\cdot \eta-   <A\eta,\eta>} d \eta= \frac{\pi^{\frac{n}{2}}}{\sqrt{\det A}} e^{-\frac{1}{4} <A^{-1} z,z>},
\end{equation} and
\begin{equation}\label{intA-1}
\int_{\mathbb R^n} e^{z\cdot \eta-   <A^{-1}\eta,\eta>} d \eta= \pi^{\frac{n}{2}} \sqrt{\det A} e^{\frac{1}{4} <Az,z>},
\end{equation}
 for any $z\in \mathbb R^n$ where $<A\eta,\eta>$ stands for $\eta^T A\eta$ for any $\eta\in\mathbb R^n$.
  \end{lemma}

\noindent  {\bf Proof of Lemma \ref{propF}}.  Note that $\div(A\nabla u) = \sum_{i,j=1}^n a_{ij} u_{x_ix_j}$. Let the notation $\mathcal F$ stand for the Fourier transform that is
\begin{equation*}\label{}
 \mathcal F(u)(\zeta,t)=(2\pi)^{-\frac{n}{2}} \int_{\mathbb R^n} e^{-i x\cdot \zeta} u(x,t) dx.
 \end{equation*}
From properties of Fourier transform we have
\begin{equation*}\label{}
\mathcal F( u_{x_ix_j})(\zeta,t)= - \zeta_i \zeta_j \mathcal F(u)(\zeta,t) \ \ \ \text{and} \ \ \ \mathcal{F} ( q\cdot \nabla u)(\zeta,t)= i  q\cdot \zeta \mathcal F(u)(\zeta,t).
 \end{equation*}
Applying the Fourier transform to  (\ref{mainl}) we get
\begin{equation*}\label{}
\partial_t \mathcal F( u)(\zeta,t) +  i  q\cdot \zeta \mathcal F(u)(\zeta,t)=-\sum_{i,j} a_{ij}\zeta_i \zeta_j \mathcal F(u)(\zeta,t)+f'(0)   \mathcal F(u)(\zeta,t).
 \end{equation*}
This is a first order linear differential equation and the solution
is
%This is a first order differential equation with an initial value $\mathcal F(u)(\zeta,0)=\mathcal F(u(x,0))$. Applying the method of characteristics we get
\begin{equation*}\label{}
\mathcal F(u)(\zeta,t) = %g'(0) \mathcal F(M_m)(\zeta)
\mathcal F(u)(\zeta,0) e^{ tf'(0)   - t\left(  i q\cdot \zeta+\sum_{i,j} a_{ij}\zeta_i \zeta_j \right)}.
 \end{equation*}
For some $k(x,t)$, define
\begin{equation*}\label{}
\mathcal F(k)(\zeta,t):= (2\pi)^{-\frac{n}{2}}   e^{ tf'(0)  - t\left(  i q\cdot \zeta+\sum_{i,j} a_{ij}\zeta_i \zeta_j \right)}.
 \end{equation*}
Therefore,
$$\mathcal F(u)(\zeta,t)=(2\pi)^{\frac{n}{2}}
%g'(0) \mathcal F(M_m)(\zeta)
\mathcal F(u)(\zeta,0) \mathcal F(k)(\zeta,t). $$ From the properties of the Fourier transform we have
\begin{equation}\label{uxt}
u(x,t)=%g'(0)\int_{\mathbb R^n} k(x-y,t) M_m(y) dy  ,
\int_{\mathbb R^n} k(x-y,t) u(y,0) dy  ,
\end{equation}
 where
 \begin{eqnarray*}\label{}
k(x,t) &=& e^{ tf'(0)} (2\pi)^{-n}  \int_{\mathbb R^n} e^{i x\cdot \zeta}   e^{- t\left(  i q\cdot \zeta+\sum_{i,j} a_{ij}\zeta_i \zeta_j \right)} d\zeta,  \\
&=& e^{ tf'(0)} (2\pi)^{-n} \int_{\mathbb R^n} e^{i (x-t q)\cdot \zeta}   e^{- t\sum_{i,j} a_{ij}\zeta_i \zeta_j} d\zeta.
  \end{eqnarray*}
  Now set $\eta=\zeta \sqrt{t}$. Then
 \begin{equation*}\label{}
k(x,t) = e^{ tf'(0)} (2\pi)^{-n} t^{-\frac{n}{2}} \int_{\mathbb R^n} e^{i \left(\frac{x-t q}{\sqrt{t}}\right)\cdot \eta-<A\eta,\eta>} d\eta.
\end{equation*}
Applying (\ref{intA}) with $z=\frac{x-t q}{\sqrt{t}}$  for any $x\in\mathbb R^n$ and  $t\in \mathbb R^+$ we get the following explicit formula for $k$
\begin{equation*}
k(x,t) = e^{ tf'(0)} (2\pi)^{-n} t^{-\frac{n}{2}}  \frac{\pi^{\frac{n}{2}}}{\sqrt{\det A}}  e^{-\frac{1}{4} <A^{-1} \left(\frac{x-t q}{\sqrt{t}}\right),\left(\frac{x-t q}{\sqrt{t}}\right)>} .
\end{equation*}
We now substitute the above formula for $k$ in (\ref{uxt}) to obtain
\begin{equation*}
u(x,t)=
e^{ tf'(0)} (2\pi)^{-n} t^{-\frac{n}{2}}  \frac{\pi^{\frac{n}{2}}}{\sqrt{\det A}}  \int_{\mathbb R^n}  e^{-\frac{1}{4} <A^{-1} \left(\frac{x-y-t q}{\sqrt{t}}\right),\left(\frac{x-y-t q}{\sqrt{t}}\right)>} u(y,0) dy  .
\end{equation*}
We apply the above integral operator to establish the formulas for both $L$ and $\bar L$ in (\ref{ml}). Let $u^{(m)}(x,t)$ be a solution of the linear equation (\ref{mainl})  with initial value (\ref{maincl}) that is   $ u^{(m)}(x,0) = g'(0) M_m(x) $. Then,
\begin{equation*}
u^{(m)}(x,t)=g'(0) e^{ tf'(0)} (2\pi)^{-n} t^{-\frac{n}{2}}  \frac{\pi^{\frac{n}{2}}}{\sqrt{\det A}}  \int_{\mathbb R^n}  e^{-\frac{1}{4} <A^{-1} \left(\frac{x-y-t q}{\sqrt{t}}\right),\left(\frac{x-y-t q}{\sqrt{t}}\right)>} M_m(y)dy.
\end{equation*}
From (\ref{maincl2}) that is $M_{m+1}(x)=u^{(m)}(x,1)$ by setting  $t=1$ in the above we conclude
\begin{eqnarray*}
\nonumber M_{m+1}(x) &=& g'(0) e^{ f'(0)} (2\pi)^{-n}  \frac{\pi^{\frac{n}{2}}}{\sqrt{\det A}}  \int_{\mathbb R^n}  e^{-\frac{1}{4} <A^{-1} \left(x-y- q\right),\left(x-y- q \right)>} M_m(y)dy
\\&=& L[M_m(x)].
\end{eqnarray*}
Note that the operator $L$ is defined on the set of all continuous functions as
\begin{eqnarray*}
 \nonumber L(v)(x) &:=& g'(0) e^{ f'(0)} (2\pi)^{-n}  \frac{\pi^{\frac{n}{2}}}{\sqrt{\det A}}  \int_{\mathbb R^n} e^{-\frac{1}{4} <A^{-1} \left(y- q\right),\left(y- q \right)>} v(x-y) dy
\\ &=& \int_{\mathbb R^n}  v(x-y) m(y,dy),
\end{eqnarray*}
where the measure $m$ is defined as
\begin{equation*}
m(y,dy) := g'(0) e^{ f'(0)} (2\pi)^{-n}  \frac{\pi^{\frac{n}{2}}}{\sqrt{\det A}}  e^{-\frac{1}{4} <A^{-1} \left(y- q\right),\left(y- q \right)>}.
\end{equation*}
Now, let $u^{(m)}(x,t)$ be a solution of the linear equation (\ref{mainl})  with initial value (\ref{mainbl}) that is   $ u^{(m)}(x,0) = \bar M_{m}(x) $. Then,
\begin{equation*}
u^{(m)}(x,t)= e^{ tf'(0)} (2\pi)^{-n} t^{-\frac{n}{2}}  \frac{\pi^{\frac{n}{2}}}{\sqrt{\det A}}  \int_{\mathbb R^n}  e^{-\frac{1}{4} <A^{-1} \left(\frac{x-y-t q}{\sqrt{t}}\right),\left(\frac{x-y-t q}{\sqrt{t}}\right)>} \bar M_m(y)dy.
\end{equation*}
From (\ref{maincl2}) that is
$$\bar M_{m+1}(x)=g'(0)  \int_{\mathbb R^n} K(x-z)u^{(m)}(z,1) dz , $$ by setting  $t=1$ in the above we conclude
\begin{eqnarray*}
\nonumber \bar M_{m+1}(x) &=& g'(0) e^{ f'(0)} (2\pi)^{-n}  \frac{\pi^{\frac{n}{2}}}{\sqrt{\det A}}  \int_{\mathbb R^n}  \int_{\mathbb R^n}  K(x-z) e^{-\frac{1}{4} <A^{-1} \left(z-y- q\right),\left(z-y- q \right)>} M_m(y)dy dz
\\&=&\nonumber \bar L[\bar M_m(x)].
\end{eqnarray*}
Note that the operator $\bar L$ is defined on the set of all continuous functions as
\begin{eqnarray*}
 \nonumber \bar L(v)(x) &:=& g'(0) e^{ f'(0)} (2\pi)^{-n}  \frac{\pi^{\frac{n}{2}}}{\sqrt{\det A}} \int_{\mathbb R^n}  \int_{\mathbb R^n} K(z-y)  e^{-\frac{1}{4} <A^{-1} \left(z- q\right),\left(z- q \right)>}  v(x-y) dz dy
\\ &=& \int_{\mathbb R^n}  v(x-y) l(y,dy),
\end{eqnarray*}
where the measure $l$ is defined as
\begin{equation*}
l(y,dy) := g'(0) e^{ f'(0)} (2\pi)^{-n}  \frac{\pi^{\frac{n}{2}}}{\sqrt{\det A}} \int_{\mathbb R^n} K(z-y)  e^{-\frac{1}{4} <A^{-1} \left(z- q\right),\left(z- q \right)>} dz.
\end{equation*}
This completes the proof.

\hfill $ \Box$
 \\
 \\
 \noindent {\bf Proof of Lemma \ref{lemfor}}.  We now compute the integral in the right-hand side of (\ref{measm})
 \begin{eqnarray*}\label{}
\nonumber \int_{\mathbb R^n} e^{s x\cdot e}  m(x,dx)  &=& g'(0) e^{ f'(0)} (2\pi)^{-n}  \frac{\pi^{\frac{n}{2}}}{\sqrt{\det A}}  \int_{\mathbb R^n} e^{s x\cdot e -\frac{1}{4} <A^{-1} \left(x- q\right),\left(x- q \right)>} dx
\nonumber
\\&=&  g'(0) e^{ f'(0)} (2\pi)^{-n}  (2)^n \frac{\pi^{\frac{n}{2}}}{\sqrt{\det A}} \int_{\mathbb R^n} e^{s e\cdot (2 \eta + q) -<A^{-1} \eta,\eta>} d\eta
\nonumber \\&=& g'(0) e^{ f'(0)} (\pi)^{-n} e^{\mu e \cdot  q}  \frac{\pi^{\frac{n}{2}}}{\sqrt{\det A}} \int_{\mathbb R^n} e^{\eta\cdot(2s e)  -<A^{-1} \eta,\eta>} d\eta .
  \end{eqnarray*}
Applying  formula (\ref{intA-1})  for $z=2\mu e$ we obtain
\begin{eqnarray*}\label{formul}
\nonumber \int_{\mathbb R^n} e^{s x\cdot e}  m(x,dx)  &=& g'(0) e^{ f'(0)} (\pi)^{-n} e^{s e \cdot  a}  \frac{\pi^{\frac{n}{2}}}{\sqrt{\det A}} \pi^{\frac{n}{2}} \sqrt{\det A} e^{\frac{1}{4} <A(2s e),(2s e)>}
\\&=& g'(0) e^{ f'(0)} e^{s e \cdot  q+s^2 <Ae,e>}.
  \end{eqnarray*}
    We now compute the integral in the right-hand side of (\ref{measl})
 \begin{eqnarray*}\label{}
\nonumber \int_{\mathbb R^n} e^{s x\cdot e}  l(x,dx)  &=& g'(0) e^{ f'(0)} (2\pi)^{-n}  \frac{\pi^{\frac{n}{2}}}{\sqrt{\det A}}  \int_{\mathbb R^n}\left[ \int_{\mathbb R^n} e^{s x\cdot e} K(y-x) dx \right]e^{ -\frac{1}{4} <A^{-1} \left(y- q\right),\left(y- q \right)>} dy
\\&=&
 g'(0) e^{ f'(0)} (2\pi)^{-n}  \frac{\pi^{\frac{n}{2}}}{\sqrt{\det A}}  k(s)\int_{\mathbb R^n} e^{s y\cdot e -\frac{1}{4} <A^{-1} \left(y- q\right),\left(y- q \right)>} dy  ,
  \end{eqnarray*}
where $k(s):=\int_{\mathbb R^n} e^{sz\cdot e} K(-z) dz$. Applying the same arguments at in the above for the measure $m$ completes the proof.

\hfill $ \Box$
\\
\\
\noindent  {\bf Proof of Theorem \ref{thmspk}}. Let $K(x)$ be the Gaussian distribution with mean $\mu$ and covariance matrix $B$ denoted by
\begin{equation*}\label{}
K_G(x) : =\frac{1}{(4\pi)^{n/2} \sqrt{\det B}} e^{-\frac{1}{4} <B^{-1}(x-\mu),(x-\mu)>}.
\end{equation*}
Here,  $\mu$ is a real vector and $B$ stands for some positive definite matrix.   In order to establish the formula of the traveling wave speed $c^*(e)$ for the linear operator $\bar L$, we minimize the following function for $s>0$, as it is given in Lemma \ref{thmsp1} that is
\begin{equation*}\label{}
 c^*(e) = \inf_{s>0} \left\{ \frac{1}{s} \ln \int_{\mathbb R^n} e^{s x\cdot e} l(x,dx) \right\},
 \end{equation*}
where the measure $l$ is given by Lemma \ref{propF},
\begin{equation*}\label{}
l(y,dy) = g'(0) e^{ f'(0)} (2\pi)^{-n}  \frac{\pi^{\frac{n}{2}}}{\sqrt{\det A}} \int_{\mathbb R^n} K(x-y)  e^{-\frac{1}{4} <A^{-1} \left(x- q\right),\left(x- q \right)>} dx.
 \end{equation*}
Define $W_G(s):=\frac{1}{s} \ln \int_{\mathbb R^n} e^{s x\cdot e} l(x,dx)$ for $s>0$. Note that from  Lemma  \ref{lemfor}, we conclude that
$$ W_G(s)= \frac{\ln k(s)}{s}+ \frac{\ln[g'(0) e^{ f'(0)} ]}{s} + e \cdot  q+s <Ae,e>, $$
where $k(s)=\int_{\mathbb R^n} K_G(x) e^{-s x\cdot e} dx$. We now compute the latter integral to evaluate $k(s)$, that is
\begin{equation*}\label{}
 k(s)
=  \frac{1}{(4\pi)^{n/2} \sqrt{\det B}} \int_{\mathbb R^n} e^{-s x\cdot e}  e^{-\frac{1}{4} <B^{-1}(x-\mu),(x-\mu)>} dx. 
  \end{equation*}
We now set a change of variable $\eta=\frac{x-\mu}{2}$. Therefore, Lemma \ref{lemcom} implies
\begin{eqnarray*}\label{}
 k(s) &=&  \frac{1}{\pi^{n/2} \sqrt{\det B}} e^{-s\mu\cdot e} \int_{\mathbb R^n} e^{(-2se)\cdot \eta- <B^{-1}\eta,\eta>}  d\eta
 \\&=&   e^{-s\mu\cdot e} e^{\frac{1}{4} <B(-2se),(-2se)>}
  \\&=&  e^{-s\mu\cdot e+<Be,e>s^2}.
  \end{eqnarray*}
  We now rewrite $W_G(s)$ as
  $$ W_G(s)=   \frac{\ln[g'(0) e^{ f'(0)} ]}{s} + e \cdot ( q-\mu) +s <(A+B)e,e>.$$
 Therefore,
\begin{equation*}\label{}
 c_{K_G}^*(e) = \inf_{s>0} W_G(s).
 \end{equation*}
  It is straightforward to compute the minimizer of the function $W_G$ as
 $$s^*:=\sqrt{\frac{\ln[g'(0)]+ f'(0) }{<(A+B)e,e>} }.$$
Therefore, $c_{K_G}^*(e)=W_G(s^*)$ and this completes the proof.

\hfill $ \Box$
\\
\\
%\noindent  {\bf Proof of Lemma \ref{lempers}}. We only provide the proof for the integral operator $\bar L$. In  the proof of Proposition \ref{propF}, we arrived at
%\begin{equation} u^{(m)}(x,t)=%g'(0) e^{ tf'(0)} (2\pi)^{-n} t^{-\frac{n}{2}}  \frac{\pi^{\frac{n}{2}}}{\sqrt{\det A}}  \int_{\mathbb R^n}  e^{-\frac{1}{4} <A^{-1} \left(\frac{x-y-t q}{\sqrt{t}}\right),\left(\frac{x-y-t q}{\sqrt{t}}\right)>} M_m(y)dy e^{ tf'(0)} (2\pi)^{-n} t^{-\frac{n}{2}}  \frac{\pi^{\frac{n}{2}}}{\sqrt{\det A}}  \int_{\mathbb R^n}  e^{-\frac{1}{4} <A^{-1} \left(\frac{x-y-t q}{\sqrt{t}}\right),\left(\frac{x-y-t q}{\sqrt{t}}\right)>} u^{(m)}(y,0) dy  .\end{equation}\hfill $ \Box$\\ \\
\noindent  {\bf Proof of Theorem \ref{thmcw}}. The operators $P$ and $Q$ and the measures  $m$ and $l$  satisfy all of the assumptions of Theorem 6.6 provided by  Weinberger in  \cite{w82}  in regards to the existence of traveling wave. More precisely,  as it was discussed in Section \ref{secspeedn}, the operators $P$ and $Q$ satisfy the assumptions (H1)-(H5).  From the assumption (G0) on the function $g$,  we have  $\frac{g(s)}{s}$ is a nonincreasing function of $s$ in $\mathbb R^+$ that is
  \begin{equation*}
\frac{g(\alpha)}{\alpha} \ge \frac{g(s)}{s} \ \ \ \ \text{for}  \ \ \ \ s\ge \alpha.
\end{equation*}
  Now set $\alpha=\rho s \le s$ for $0\le \rho \le 1$ and apply the above inequality to get
  \begin{equation*}
g(\rho s) \ge \rho g(s) .
\end{equation*}
This implies that for all $0\le \rho \le 1$ we have
\begin{equation*}
Q[\rho u] \ge \rho Q[u] ,
\end{equation*}
where $u$ is a solution of (\ref{main}) and $Q$ is given by (\ref{op}). In the light of arguments in Section \ref{secprop}, we choose a function $\phi$  that satisfies the corresponding properties.  For each positive integer $k$ we define the sequence
$$ a_{m+1}(c,e;k,s) : =\max\{k^{-1}\phi(s), Q[a_m(c,e,k;x\cdot e+s+c)] (0)\},$$
when $a_0(c,e,k;s)=k^{-1}\phi(s)$.  The sequence $a_n$ is nonincreasing in $c$, $k$, and $s$ and nondecreasing in $n$. As $n\to\infty$, it converges to a limit function $a$ which is nonincreasing in $c$, $k$ and $s$. The function  $a(c,e,k;x\cdot e+t)$ is continuous in terms of $x$ and
 $\lim_{s\to-\infty} a(c,e,k;s)=\pi_1$ and $\lim_{s\to\infty} a(c,e,k;s)=0$ for $c\ge c^*(e)$. Furthermore,
\begin{equation}\label{acek}
 a(c,e;k,s) : =\max\{k^{-1}\phi(s), Q[a(c,e,k;x\cdot e+s+c)] (0)\} .
 \end{equation}
Set $x_0\in\mathbb R^n$ such that $x_0\cdot e>0$. For every integer $j$ and any $c\ge c^*(e)$ define the sequence
$$ T_k(l)=\frac{1}{2}\left[a(c,e,k;lx_0\cdot e)+a(c,e,k;(l+1)x_0\cdot e)  \right] . $$
Then, $T_k(l)$ is nonincreasing in $l$,  $T_k(-\infty)=\pi_1$ and $T_k(\infty)=0$.  Since $a$ decreases
from $\pi_1$ to $0$ as $s$ goes from $-\infty$ to $\infty$,
$$ T_k(l) -T_k(l-1) =\frac{1}{2}\left[a(c,e,k;(l+1)x_0\cdot e)-a(c,e,k;(l-1)x_0\cdot e)  \right] \le \frac{\pi_1}{2}  . $$
Therefore,  there is an integer $l_k$ such that
$$ \frac{\pi_1}{4} \le T_k(l_k) \le \frac{3\pi_1}{4}. $$
 We now consider the sequence $a(c,e,k;x\cdot e+l_k x_0\cdot e)$ for $k\in\mathbb  N$. There is a subsequence $k_i$ of the integers such that $a(c,e,k_i; x\cdot e+l_{k_i}x_0 \cdot e)$ converges uniformly for $x$ on
bounded subsets of $\mathbb R^n$ to a function $W(x\cdot e)$ defined on $\mathbb R^n$.  From this we conclude that there exists a
subsequence $k^{(1)}_i$ such that $a(c,e,k_i^{(1)}; x\cdot e +l_{k_i^{(1)}} x_0\cdot e +c)$ also converges uniformly on
bounded subsets of $\mathbb R^n$ to a function $W(x\cdot e+c)$.  Applying this argument we establish a sequence $a(c,e,k_i^{(m)}; x\cdot e +l_{k_i^{(m)}} x_0\cdot e +mc)$ that is uniformly convergent to  a function $W(x\cdot e+mc)$ for any $m\in\mathbb N$. Using a diagonal process we arrive at a sequence $\bar k_i$ and  taking limits in (\ref{acek})  for $k=\bar k_i$ and $s=y\cdot e +l_{\bar k_i} x_0 \cdot e -(m+1)c  $ we get
 $$ W(y\cdot e -(m+1)c)=Q[W(x\cdot e - mc)] (y). $$
 Therefore, $u_m(x)=W(x\cdot e-mc)$ is a traveling wave solution of $u_{m+1}=Q[u_m]$.

\hfill $ \Box$

\noindent {\bf Acknowledgement:}  The authors would like to thank the anonymous referees for many valuable suggestions.

 %%%%%%%%%%%%%%%%%%%%%%%%

\end{document}